\documentclass[11pt]{amsart}
\usepackage{amsmath}
\usepackage{amssymb}
\usepackage{latexsym}
\usepackage{a4wide}
\usepackage{epsfig}

\input{psfig}

\newcommand{\N}{\ensuremath{\mathbb N}}
\newcommand{\Z}{\ensuremath{\mathbb Z}}
\newcommand{\R}{\ensuremath{\mathbb R}}

\newcommand{\C}{\ensuremath{\mathbb C}}
\newcommand{\Q}{\ensuremath{\mathbb Q}}

\newcommand{\X}{\ensuremath{\mathbb X}}

\newcommand{\D}{\ensuremath{{\mathcal D}}}

\newcommand{\T}{\ensuremath{{\mathcal T}}}
\newcommand{\A}{\ensuremath{{\mathcal A}}}

\newcommand{\Xc}{\ensuremath{{\mathcal X}}}

\renewcommand{\rho}{\varrho}
\renewcommand{\phi}{\varphi}

\DeclareMathOperator{\inn}{int}

\DeclareMathOperator{\cl}{cl}

\newtheorem{thm}{Theorem}[section]
\newtheorem{defi}[thm]{Definition}

\newtheorem{cor}[thm]{Corollary}
\newtheorem{lem}[thm]{Lemma}

\begin{document}

\title{Self-dual tilings with respect to star-duality}

\author{D. Frettl\"oh}
\address{Fakult\"at f\"ur Mathematik, Universit\"at Bielefeld, 
Postfach 100131, 33501 Bielefeld,  Germany}
\email{dirk.frettloeh@math.uni-bielefeld.de}
\urladdr{http://www.math.uni-bielefeld.de/baake/frettloe}

\begin{abstract} 
The concept of star-duality is described for self-similar
cut-and-project tilings in arbitrary dimensions. This generalises 
Thurston's concept of a Galois-dual tiling. The dual tilings of the
Penrose tilings as well as the Ammann-Beenker tilings are calculated.
Conditions for a tiling to be self-dual are obtained. 
\end{abstract} 

\maketitle

\section{Introduction}
Nonperiodic discrete structures like infinite words, tilings, or point
sets, are a rich source of exciting objects studied
within many fields of mathematics, including automata theory,
combinatorics, discrete geometry, dynamical systems, group theory,
mathematical physics and number theory. Many interesting examples
of such structures can be produced by substitutions as well as
by cut-and-project schemes. For such structures there are several
concepts of 'dual' tilings (words, point sets) \cite{sai},
\cite{sirw}, \cite{vin}. The present article describes
how Thurston's concept of Galois-duality \cite{thu} generalises to 
arbitrary cut-and-project sets. This generalisation does not depend on
Galois conjugacy but on the star map (see \cite{moo}, or Definition
\ref{defmodset} below) wherefore the term star-dual (resp.\
$\star$-dual) is used, rather than Galois dual. The benefits of
star-duality are that it gives a purely algebraic approach to certain
geometric problems; that it covers a pretty wide range of examples;
and that it may help to friend the certain concepts of
duality. Furthermore, and hopefully, it can be  utilised to attack
open problems in the theory of non-periodic structures. 

The first section of the present paper  
collects necessary facts and definitions about substitutions and
cut-and-project tilings, resp.\ model sets. Section \ref{secstardual} 
contains the definition of $\star$-duality, together with an example. 
The calculation of the $\star$-dual tilings of the famous 
Penrose tilings and Ammann-Beenker tilings \cite{gs} is carried
out in Section \ref{secdim2}. It turns out that the star-dual
tilings of the Penrose tilings are the well-known T\"ubingen
triangle tilings (and vice versa), which coincides with a related
result  in \cite{bksz}. In contrast, the star-dual tilings
of the Ammann-Beenker tilings are very similar to the 
Ammann-Beenker tilings themselves. In fact, they are mutually 
locally derivable \cite{bsj} with the Ammann-Beenker tilings. 
This motivates the question whether there are self-dual tilings
w.r.t.\ star-duality. Therefore, Section \ref{secselfd} is dedicated
to self-duality in arbitrary dimensions. 
A necessary, and a sufficient condition for a tiling to be
self-dual  w.r.t.\ star-duality are obtained.

The concept of star-duality is not new, see for instance \cite{thu},
\cite{gel}, \cite{vin}. The relevant results of this
paper are the following: First, Theorem \ref{D+t} gives a useful
technical result which is exploited in the next sections, and will be
exploited further in future work. Second, the dual tilings of the
Penrose tilings and the Ammann-Beenker tilings are obtained, where the
latter is done for the first time here, up to knowledge of the author,
where the former one is obtained in \cite{bksz} with much more effort
in a slightly different framework (but not in \cite{gel}, as the title
might suggest).  
Last, Lemma \ref{lemsd} and Corollary \ref{corsd} are first steps
towards a classification of self-dual tilings, in one dimension as
well as in arbitrary dimension. 

{\bf Remark:} For convenience, we will exclude periodic structures
in the following. Each occurring biinfinite sequence, tiling etc.\ is
assumed to be nonperiodic. The set of positive integers is denoted by
$\N$. The closure (interior) of a set $A$ is denoted by $\cl(A)$
($\inn(A)$). 

\subsection{\bf  Substitution sequences} \label{substseq}
Let $\A$ be an {\em alphabet}, i.e., a finite set of letters; let
$\A^{\ast}$ be the set of all finite words over the alphabet $\A$,
i.e., all words which are concatenations of letters of
$\A$. Usually, the empty word $\epsilon$ (consisting of no letters) is
required to be an element of $\A^{\ast}$. A {\em substitution} is a
map $\zeta : \A \to \A^{\ast} \setminus \{ \epsilon \}$. We will call
such a substitution also {\em symbolic substitution} if we want to
emphasise the difference to a tile-substitution. 
By setting $\zeta(a_1 \cdots a_m) = \zeta(a_1) \cdots \zeta(a_m)$,
$\zeta$ acts also as a map from $\A^{\ast}$ to $\A^{\ast}$, from
$\A^{\N}$ to $\A^{\N}$, and from $\A^{\Z}$ to $\A^{\Z}$. 
The family of all biinfinite words,
where each finite subword is contained in some $\sigma^k(a)$ 
($a \in \A$), is denoted $\Xc_{\sigma}$. There is a
lot of literature devoted to symbolic substitutions, see for instance
\cite{fogg} and references therein for more details.

\subsection{\bf Substitution tilings} \label{substtil}

A tile is a nonempty compact subset of $\R^d$ which is the closure of
its interior. A tiling of $\R^d$ is a set 
$\{ T^{}_j \}^{}_{j \in \N}$ of tiles  which is a covering of
$\R^d$ (i.e., $\R^d = \bigcup_{j \in \N} T^{}_j$) that is
non-overlapping (i.e., $\inn(T^{}_j) \cap \inn(T^{}_i) =
\varnothing$ for $j \ne i$). A simple way to generate tilings ---
periodic ones as well as nonperiodic ones --- is via a
tile-substitution.  Roughly speaking, a tile-substitution is given by
a finite set  $\{T^{}_1, \ldots, T^{}_m \}$ of tiles, the 
{\em  prototiles}, an expanding map 
$Q$, and a rule how to dissect each $QT^{}_j$ into isometric copies of
some prototiles  $T^{}_i$. 

\begin{figure}
\epsfig{file=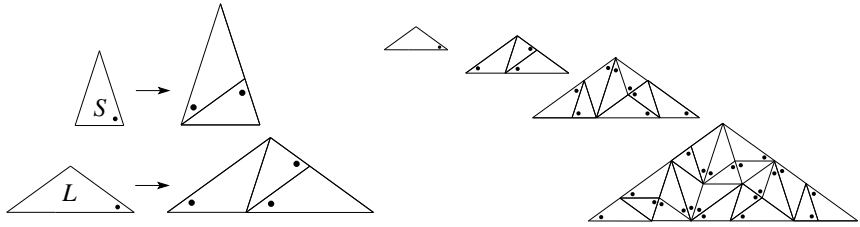} 
\caption{ \label{fig:penr-subst}
The substitution for the Penrose tiling in the version
  with triangles $S$ and $L$ (left), the first iterates of the
  substitution on the obtuse triangle $L$ (right). The substitution
  factor in this example is $\tau = \frac{\sqrt{5}+1}{2}$.}
\end{figure}

{\bf Example:} A diagram of the tile-substitution for the Penrose
tiling \cite{gs} is shown in Figure \ref{fig:penr-subst}. Here, we use
triangles as prototiles rather than rhombi or the celebrated kites and 
darts. Since the two 
prototiles are mirror-symmetric, but their substitutions are not, we
have to respect the orientation of each tile. This is achieved by a
black dot on each triangle. As indicated in
the figure, iterating the substitution several times fills larger and
larger portions of space, eventually resulting in a tiling of the plane.
However, the proper definition of a
tile-substitution and a substitution tiling is as follows.

\begin{defi} \label{defsubst}
A (self-affine) tile-substitution is defined via a set of
prototiles and a map $\sigma$. Let $T^{}_1, T^{}_2, \ldots
T^{}_m$ be nonempty compact 
sets --- the {\em prototiles} --- in $\R^d$, such that 
$\cl( \inn(T^{}_j) ) =  T^{}_j$ for each $i \le m$. Let $Q$ be an
expanding linear map such that 
\begin{equation} \label{eq:eifs}
 Q T^{}_j = \bigcup_{i=1}^m T^{}_i + \D^{}_{ij}  \quad (1 \le j \le m), 
\end{equation}
where the union is not overlapping (i.e., the interiors of the tiles
in the union are pairwise disjoint); and each $\D^{}_{ij}$ is a finite
(possibly empty) subset of $\R^d$, called {\em digit set}. Then  
\[ \sigma( T^{}_j ) := \{ T^{}_i + \D^{}_{ij} \, | \, i=1 \ldots m
\} \]
is called a {\em (tile-)substitution}.  
By $\sigma(T^{}_j+x):=\sigma(T^{}_j)+Q x$ and $\sigma(\{T,T'\}) := \{
\sigma(T), \sigma(T') \}$, $\sigma$ extends in a natural way to all
sets $\{ T^{}_{i(k)} + x^{}_k \, | \, k \in I \}$.   
\end{defi}
The sum $T+\D^{}_{ij} = T + \{d^{}_{ij1}, \ldots, d^{}_{ijk} \}$ means
$\{ T+ d^{}_{ij1}, \ldots, T+ d^{}_{ijk} \}$, with the convention
$T + \varnothing = \varnothing$, if  $\D^{}_{ij}$ is empty. 
The set $I$ in the last
line of the definition  can be any index set. In particular, if
$I = \N$, such a set can describe a tiling, given
that the tiles do not  overlap, and the union of all tiles is the
entire space $\R^d$. Before stating the definition of a substitution
{\em tiling}, let us make two remarks.  

{\bf Remark:} Often, a tile-substitution is defined without the
requirement \eqref{eq:eifs}. To be precise, 
if the support of $\sigma(T_j)$ equals $QT_j$ --- that is, if 
\eqref{eq:eifs} holds ---, the 
substitution is called {\em self-affine}. Moreover, 
if $Q$ is a homothety (i.e., just a multiplication
by a scalar factor $\lambda$), the substitution is called 
{\em self-similar}. In this case, $\lambda$ is called the 
{\em substitution factor} of the substitution. 
E.g., the expanding linear map for the Penrose tiling is $Q = \bigl(
\begin{smallmatrix} \tau & 0 \\ 0 & \tau \end{smallmatrix} \bigr)$,
where $\tau=(\sqrt{5}+1)/2$ is the golden mean. Thus,
the substitution factor of the Penrose substitution in
Figure \ref{fig:penr-subst} is $\tau$. 
Whenever the term 'self-similar tiling' ('self-affine tiling') is 
used in the sequel this is to be understood as 'self-similar
substitution tiling' ('self-affine substitution tiling').

{\bf Remark:} Every symbolic substitution (see Section \ref{substseq})
gives rise to a self-similar tile-substitution in $\R^1$ in a canonical way,
just by replacing each letter by an interval. This can be done in such a
way that the tile-substitution is self-similar. Namely, one uses
the entries of a positive PF-eigenvector of the substitution matrix,
see \cite{fogg}, \cite{fre}. 

\begin{defi} \label{substtildef}
Let $\sigma$ be a tile-substitution.
A tiling $\T$ is called a {\em substitution tiling}
(with substitution $\sigma$), if each finite subset of $\T$ is a
translate of some subset of some $\sigma^k(T^{}_j)$, where $k \ge 0$,
and $T^{}_j$ some prototile.  The family of
all substitution tilings with substitution $\sigma$ is the 
{\em  tiling space} $\X^{}_{\sigma}$.  
\end{defi}
This definition of a tiling space is equivalent
to the usual one \cite{sol}, whenever $\sigma$ is {\em primitive}. 
Primitive means that there is $k \ge 0$ such that a power of the 
{\em substitution matrix} (or 'incidence matrix') 
$M_{\sigma}=(|\D^{}_{ij}|)^{}_{1 \le i,j \le m}$ is strictly positive.
Here, $|\D^{}_{ij}|$ denotes the cardinality of $\D^{}_{ij}$.   
It is also consistent with the usual definition of the associated
dynamical system $(\Xc_{\sigma},S)$ arising from $\sigma$, see
\cite{fogg}, as well as with our definition of $\Xc_{\sigma}$
in Section \ref{substseq}.

\subsection{Iterated function systems}

The prototiles of a self-affine substitution can be derived 
from $Q$ and $(\D^{}_{ij})^{}_{1 \le i,j \le m}$ uniquely as follows: 
Multiplying the entire equation system \eqref{eq:eifs} by $Q^{-1}$ 
yields the corresponding iterated function system (IFS) for the tiles 
of $\sigma$: 
\begin{equation} \label{eq:ifs}
T^{}_j = \bigcup_{i=1}^m Q^{-1} (T^{}_i + \D^{}_{ij})  \quad 
(1 \le j \le m)
\end{equation}

Such an IFS is known to possess a unique nonempty compact solution
(see for instance \cite{ban} and references therein).  
Since the prototiles of $\sigma$ fulfil \eqref{eq:ifs}, they are the
unique solution of \eqref{eq:ifs}. In particular, by Definition
\ref{defsubst}, the IFS \eqref{eq:ifs} is {\em non-overlapping}, that
means, for each $j$, the interiors of the sets in the union are
pairwise disjoint. Note that two parts $T_i, T_j$ of the solution are
allowed to overlap. The point is that \eqref{eq:ifs} induces a
partition of each $T_j$ up to a set of measure 0.

Originally, the term IFS was coined for the case $m=1$ only. An IFS
with $m>1$ components is often called 'graph-iterated function system'
or 'multi-component IFS'. Here, we use the term IFS for all of them.

{\bf Remark:} We just mentioned that each self-affine substitution gives
rise to a non-overlapping IFS \eqref{eq:ifs}. Conversely, by multiplying
\eqref{eq:ifs} with $Q$, each non-overlapping IFS with a common
contracting map $Q^{-1}$ defines a self-affine substitution as in
Definition \ref{defsubst}. 

{\bf Remark:} An IFS as above can be written as a finite state
automaton (FSA). The states are the $T^{}_j$. The input alphabet consists
of all elements $x$ in  $\bigcup_{1 \le i,j \le m} \D^{}_{ij}$. The input
of a digit $x$  changes the state from $T^{}_j$ to $T^{}_i$, if $x \in
\D^{}_{ij}$, otherwise the state is changed to FAIL. For an example,
see Figure \ref{fig:penrauto}. For details, cf.\ \cite{thu}, \cite{gel}.   

The following theorem ensures that it does not matter where the
prototiles for a substitution $\sigma$ are located in space. The
tiling space defined by $\sigma$ depends only on the way {\em how}
tiles are dissected  into smaller pieces, not {\em where} they are
located. In particular, it  does not matter whether the prototiles
overlap or not. The important property is that the dissection of each
prototile is non-overlapping. This theorem will be useful in the next
sections, as well as for future work. 

\begin{thm} \label{D+t}
Let $\D$ be the digit set matrix for an IFS in $\R^d$ with a common
contraction $Q$. If $(A^{}_1,A^{}_2, \ldots, A^{}_m)$ is a solution
for $\D$, then $(A^{}_1 + t^{}_1, A^{}_2 + t^{}_2, \ldots, A^{}_m +
t^{}_m)$ is a solution for 
$\{\D_{ij}+(t^{}_j-Qt^{}_i)\}^{}_{1 \le i,j \le m}$ for
all $(t^{}_1, t^{}_2, \ldots, t^{}_m) \in (\R^d)^m$. 
\end{thm}
\begin{proof}
Throughout the proof, let $(A^{}_1,\ldots,A^{}_m)$ be a solution for
$\D$, i.e., 
\[ A^{}_j = \bigcup_{i=1}^m Q^{-1} (A^{}_i + \D^{}_{ij})  \quad 
(1 \le j \le m). \]

First, let  $A'_j = A^{}_j + t$ for some common $t$ ($1 \le j \le m$). 
Then the approach 
\[ \begin{array}{rclcccl}
A'_1 & = & Q A'_1 + \D^{}_{11} + x  & \cup & \cdots & \cup & Q A'_m +
\D^{}_{m1} + x \\
\vdots & & \vdots & & & & \vdots  \\
A'_m & = & Q A'_1 + \D^{}_{1m} + x  & \cup & \cdots & \cup & Q A'_m +
\D^{}_{mm} + x \\
\end{array} \]  
yields 
\[ A^{}_1 + t = Q A^{}_1 +Qt + \D^{}_{11} +x \cup \cdots \cup Q A^{}_m +Qt +
\D^{}_{m1} +x, \]
and it follows $x=t-Qt$. It is straightforward to check that 
$(A'_1, \ldots, A'_m) = (A^{}_1 + t, \ldots, A^{}_m +t)$ is a solution 
for $\{\D_{ij}+(t-Qt)\}^{}_{1 \le i,j \le m}$. 
  
Next, let $t^{}_1=t^{}_2 = \cdots = t^{}_k = 0, \; t^{}_{k+1}= \cdots
= t^{}_m = t$ 
and $A'_j = A^{}_j + t^{}_j$ for $k+1 \le j \le m$. Then the approach
\begin{equation} \label{eq:null+t}
\begin{array}{rclcl}
A^{}_1 & = & \bigcup_{i=1}^k Q A^{}_i + \D^{}_{i1}  & \cup &
\bigcup_{i=j+1}^m Q A'_{i} + \D^{}_{i,1} + x \\  
\vdots & & \vdots & & \vdots \\
A^{}_k & = & \bigcup_{i=1}^k Q A^{}_i + \D^{}_{ik}  & \cup &
\bigcup_{i=k+1}^m Q A'_{i} + \D^{}_{i,k} + x \\  
A'_{k+1} & = & \bigcup_{i=1}^k Q A^{}_i + \D^{}_{i,k+1} +y & \cup &
\bigcup_{i=k+1}^m Q A'_{i} + \D^{}_{i,k+1} + z \\  
\vdots & & \vdots & & \vdots \\
A'_m & = & \bigcup_{i=1}^k Q A^{}_i + \D^{}_{i,m} +y & \cup &
\bigcup_{i=k+1}^m Q A'_{i} + \D^{}_{i,m} + z \\  
\end{array} 
\end{equation}
yields $y=t, x=-Qt$ and $z=t-Qt$. With these values, it is 
straightforward to check that $(A^{}_1, \ldots,
A^{}_k, A^{}_{k+1} + t, \ldots A^{}_m+t)$ is a solution of
\eqref{eq:null+t}.

By combining these two results, we obtain the digit matrix for
$(A^{}_1+t^{}_1, A^{}_2 +t^{}_2, \ldots, A^{}_m +t^{}_m)$ successively as
follows.

$\D^{(1)} = \D+(t^{}_1-Qt^{}_1) = \{\D_{ij}+(t^{}_1-Qt^{}_1)\}^{}_{1
  \le i,j \le m}$ has the solution $(A^{}_1+t^{}_1, A^{}_2 +t^{}_1,
\ldots, A^{}_m +t^{}_1)$.\\[1mm] 
\[ \D^{(2)} = \D+(t^{}_1-Qt^{}_1) + \left( 
\begin{smallmatrix} 
0 & t^{}_2-t^{}_1 & \cdots & t^{}_2-t^{}_1 \\ 
-Q(t^{}_2-t^{}_1) & (E-Q)(t^{}_2-t^{}_1) & \cdots & (E-Q)(t^{}_2-t^{}_1) \\  
\vdots & \vdots & & \vdots \\
-Q(t^{}_2-t^{}_1) & (E-Q)(t^{}_2-t^{}_1) & \cdots & (E-Q)(t^{}_2-t^{}_1) 
\end{smallmatrix}
\right) \]  has the solution $(A^{}_1+t^{}_1, A^{}_2 +t^{}_1+(t^{}_2-t^{}_1),
\ldots, A^{}_m +t^{}_1 + (t^{}_2-t^{}_1))$ (where we use
\eqref{eq:null+t} with $y=t=t^{}_2-t^{}_1$, and $E$ denotes the
identity matrix), and so on, until
\[ \D^{(m)} = \D^{(m-1)} + \left( 
\begin{smallmatrix} 
0 & \cdots & 0 & t^{}_m-t^{}_{m-1} \\ 
\vdots & & \vdots & \vdots\\
0 & \cdots & 0 & t^{}_m-t^{}_{m-1} \\ 
-Q(t^{}_m-t^{}_{m-1}) & \cdots & -Q(t^{}_{m}-t^{}_{m-1}) & (E-Q)(t^{}_m-t^{}_{m-1}) 
\end{smallmatrix} \right) \]
has the solution $(A^{}_1+t^{}_1, A^{}_2 +t^{}_2, \ldots, A^{}_{m-1}+t^{}_{m-1},  
A^{}_m +t^{}_{m-1} + (t^{}_m-t^{}_{m-1}))$.\\
Now, $\D^{(m)}$ is of the desired form, which proves the claim. 
\end{proof}

Regarding the action of a common regular linear map $T$ on the solution 
$(A_1, \ldots, A_m)$ with respect to $\D$ rather than the action of
different translations is much simpler. The same is true if all
translation vectors $t_i$ in the previous theorem had been equal. 

\begin{lem} \label{lmalD}
Let $\D$ be the digit set matrix for an IFS in $\R^d$ with a common
contraction $Q$. Let $T: \R^d \to \R^d$ be any invertible linear map. 
If $(A^{}_1,A^{}_2, \ldots, A^{}_m)$ is a solution
for $\D$, then $(T A^{}_1, T A^{}_2, \ldots, T A^{}_m )$ 
is a solution for $T (\D)$, where $T (\D) = \big( T(\D^{}_{ij})
\big)^{}_{1 \le j,i \le m}$, $T(\D^{}_{ij}) = \{ T d^{}_{ijk} \, | \,
d^{}_{ijk} \in \D^{}_{ij} \}$.  
\end{lem}

\begin{proof}
\[ A_j = \bigcup_{i=1}^m A_i + \D_{ij} \quad \Rightarrow \quad T A_j =  
T \Big( \bigcup_{i=1}^m A_i + \D_{ij} \Big) 
=  \bigcup_{i=1}^m T A_i + T( \D_{ij}).  \]
\end{proof}

\subsection{Model sets}
In the sequel we consider self-similar tilings which arise from
{\em model sets}. Such tilings are known as cut-and-project tilings
\cite{baa}, \cite{mey}, \cite{moo}; see also \cite{fogg}, where the 
cut-and-project method for one-dimensional tilings is called 
'geometric realization'.

\begin{defi} \label{defmodset}
A model set is defined via a {\em cut and project scheme}, i.e., a
collection of spaces and mappings as follows.
  
Let $G,H$ be locally compact Abelian groups, $\Lambda$ be a lattice in 
$G \times H$ (that is, $\Lambda$ is a cocompact discrete subgroup of
$G \times H$),  $\pi_1: G \times H \to G, \; \pi_2: G
\times H \to H$ be projections, such that
$\pi_1 |_{\Lambda}$ is injective, and $\pi_2(\Lambda)$ is dense in H. 
Let $W \subset H$ be a compact set --- the {\em window} --- such that
the closure of $\inn(W)$ equals $W$.  
\[ \begin{array}{ccc} 
G & \stackrel{\pi_1}{\longleftarrow} G \times H
  \stackrel{\pi_2}{\longrightarrow} & H \\
\cup & \cup & \cup \\
V & \Lambda & W
\end{array} \]
Then 
\[ V := \{ \pi_1(x) \, | \, x \in \Lambda, \; \pi_2(x) \in W \} \]
is called a {\em model set}.

If $\mu(\partial(W))=0$, then $V$ is called {\em regular} model set. 
 
The {\em star map} is the map $\star: \pi_1(\Lambda) \to H, \;
x^{\star} = \pi_2(\pi^{-1}_1(x))$. 
\end{defi}

We should mention that in many cases in the literature holds $G= \R^d$
and $H =\R^e$. One may assume this reading this article
without loosing much information. However, the notion of star-duality 
is tailored to the general case, and this is where we leave the concept 
of Galois-duality. There are cut-and-project tilings (in fact, a lot
of them) of the line or the plane, where $H$ is the field $\Q^{}_p$ 
of $p$-adic integers, or a product of those with some Euclidean space
$\R^e$, see \cite{bms}, \cite{sie}, \cite{sing}, \cite{sing2}.

In the context of symbolic substitutions \cite{fogg}, the setting is 
very much the same, but unfortunately the terminology differs. For 
instance, the cut-and-project scheme is called geometric realization, 
the window is usually called (generalized) Rauzy fractal etc.

A model set is, by definition, a point set. There are several ways to 
assign a tiling to it. For instance, one may take the Voronoi cells
of the point set as tiles. Whenever we can obtain a tiling $\T$ from
a model set in such a way we call $\T$ a {\em cut-and-project tiling}.
To be precise: If it is possible to assign a tiling $\T$ to a 
model set $V$ such that $V$ and $\T$ are mutually locally derivable
(MLD) in the sense of \cite{bsj}, we call $\T$ a cut-and-project
tiling. Roughly spoken, MLD means that one can obtain $V$ from $\T$ by 
local replacement rules and vice versa.

The self-similar cut-and-project tilings are the objects we will
consider in the following. We should mention that there is a huge
number of such tilings, including the Fibonacci tilings and the
Tribonacci tilings \cite{fogg}, the Penrose tilings and the
Ammann-Beenker tilings \cite{gs}, the Chair tilings and the Sphinx
tilings \cite{bms}. A collection of several additional examples is
available online \cite{fh}.

\section{The $\star$-dual substitution} \label{secstardual}

It is known that the window of a self-similar model set is the
solution of an IFS, see for instance \cite{har},
\cite{sing}, \cite{sirw}. Since each non-overlapping IFS defines a 
tile-substitution (see the remark after \eqref{eq:ifs}),
several authors considered this 'dual substitution tiling' (or better, 
the substitution) induced by the IFS of the window.  
In \cite{thu}, Thurston uses a different construction based on
Galois conjugation in the field $\Q(\lambda)$, where $\lambda$
denotes the substitution factor of a self-similar tiling. One can
utilise Theorem \ref{D+t} to show that both approaches yield the same
tiling spaces. Thurston
called the tilings arising from this construction 'Galois duals'
of the original tilings. We proceed by stating the definition
of the $\star$-dual tiling of a self-similar cut-and-project tiling,
which  is a generalization of Galois-duality. 

In general, the groups $G$ and $H$ in Definition \ref{defmodset} can
be arbitrary locally compact Abelian groups. However, dealing with
substitutions and linear maps, we require $G$ and $H$ to be equipped 
with a vector space structure in the sequel. This is no rigid
restriction, since all examples in the literature fulfil this
condition. 

Furthermore, in this section we require the expanding linear map 
$Q: G \to G$ to be as follows: Let $T: G \times H \to G \times H$ be a
hyperbolic linear map (i.e., there is no eigenvalue $e$ of $T$ with
$|e|=1$) such that $T|_G$ is an expansion and $T|_H$ is a
contraction. Let $Q=T|_G$. Again, this is not really a restriction,
the majority of well-studied tilings fit into this framework.

\begin{defi} 
Let $\T$ be a self-similar cut-and-project tiling with 
substitution $\sigma$ and expansion $Q$. Let 
$\D=(\D^{}_{ij})_{1 \le i,j \le m}$ be a corresponding digit matrix as
in Definition \ref{defsubst} such that $\sigma$ is the substitution
induced by  $\D$. Furthermore, let $Q \pi_1(\Lambda) \subset
\pi_1(\Lambda)$, and each $d \in \D^{}_{ij}$ be in
$\pi_1(\Lambda)$, with
$\pi_1$ and $\Lambda$ as in Definition\ \ref{defmodset}. 
Then $(\D^{\star})^T$ defines a substitution $\sigma^{\star}$ with
expansion $Q':=(T|_H)^{-1}$,
which is called the {\em $\star$-dual substitution}.  
\end{defi}

Here, $\D^{\star} = (\D_{ij}^{\star})^{}_{1 \le i,j \le m}$, and 
$\D_{ij}^{\star} = \{ d^{\star} \, | \, d \in \D_{ij} \}$.

If the substitution is one-dimensional and self-similar with
substitution factor $\lambda$, where $\lambda$  is a unimodular
$PV$-number, then $Q'=\big( \text{diag} ( \lambda^{}_2, \lambda^{}_3,
\ldots, \lambda^{}_N) \big)^{-1}$, where $\lambda^{}_2, \ldots
\lambda^{}_N$  are the Galois conjugates of $\lambda$, see \cite{gel},
\cite{fre}. Therefore the name 'Galois duals' was coined in \cite{thu}.

{\bf Remark:} Note that the star-dual substitution also defines
self-similar cut-and-project tilings, 
with the roles of $G$ and $H$ in Definition\ \ref{defmodset} 
reversed. In particular, a star-dual tiling of some tiling in $G$ is
a tiling in $H$. There are tilings of the Euclidean
plane (or the line) with $p$-adic internal spaces $H$, see \cite{bms},
\cite{sing2}. For instance, the star-dual tiling of a plane tiling may be
a tiling in $\Q_p \times \Q_p$, where $\Q_p$ is the field of $p$-adic
integers, which is a locally compact Abelian group as well. 

{\bf Remark:} It follows immediately from the definition that 
$(\sigma^{\star})^{\star} = \sigma$. By Lemma \ref{lmalD}, we can
freely choose between applying the star map to $\D$ or $Q^{-1}\D$. The 
resulting substitutions are the same, up to a scaling of the
prototiles.

{\bf Example:} (Fibonacci squared) Consider the word substitution 
$\zeta$ on the alphabet $\{a,b\}$ given by 
$\zeta(a) = aab, \; \zeta(b)= ab$. The substitution matrix 
(or 'incidence matrix') is $M^{}_{\zeta}=\big( 
\begin{smallmatrix} 2 & 1 \\ 1 & 1
\end{smallmatrix} \big)$. It is well known --- and easy to see ---
that the PF-eigenvalue (the leading eigenvalue) of $M$ is the
substitution factor for the tile-substitution. Moreover, the entries
of a left eigenvector of $M$, corresponding to the PF-eigenvalue,
contains the relative lengths of the prototiles for a self-affine
tile-substitution \cite{fre}, \cite{fogg}. The leading eigenvalue in
this case is $\tau^2 = \frac{\sqrt{5}+3}{2}$, which is the substitution
factor of the desired substitution. Requiring
the prototiles to be intervals, the lengths can be chosen as
$\ell^{}_a=1, \; \ell^{}_b=\tau^{-1}$. (Here is some freedom, we may
choose any multiples of these lengths as well.) One appropriate
tile-substitution is given by  the digit set matrix
\[ \begin{pmatrix} \{0,1\} & \{0\} \\
\{ 2\} & \{ 1 \} \end{pmatrix} . \]
The prototiles $T^{}_a = [0,1]$ and $T^{}_b = [0,\tau^{-1}]$ arise as
the unique compact nonempty solution of the corresponding IFS
\begin{equation} \label{ifsfibsq}
\begin{split}
T^{}_a & = \tau^{-2} T^{}_a \cup \tau^{-2} (T^{}_a +1) \cup \tau^{-2}
(T^{}_b +2), \\  
T^{}_b & = \tau^{-2} T^{}_a \cup \tau^{-2} (T^{}_b +1). 
\end{split}
\end{equation}

\begin{figure}
\input{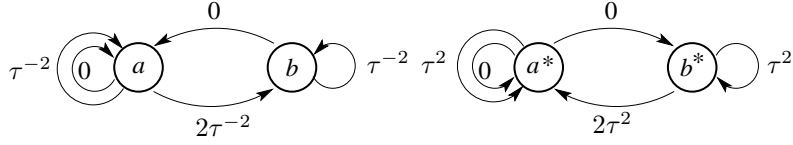}
\caption{The automaton for the IFS for the squared Fibonacci
  substitution and its dual automaton. \label{fibsqauto}}
\end{figure}

We proceed by considering the digit set of this IFS, $\D = \Big(
\begin{smallmatrix} \{0,\tau^{-2}\} & \{ 0 \} \\ \{ 2 \tau^{-2} \} &
  \{ \tau^{-2} \} \end{smallmatrix} \Big)$, in agreement with the last
remark, and in order to illustrate the action of Galois conjugation.  
The IFS \eqref{ifsfibsq} can be represented by the automaton in Figure
\ref{fibsqauto}, left. (Vice versa, the automaton gives rise to the
IFS \eqref{ifsfibsq}). The dual substitution $\zeta^{\star}$ arises
from the IFS defined by $(\D^{\star})^T$, or equivalently: from the
dual automaton, where the orientation of each edge is reversed and the
star map is applied to each edge label, see Figure \ref{fibsqauto},
right. In this case, the star map is uniquely defined by 
$1^{\star} = 1,\; \tau^{\star} = -\tau^{-1}$ (each $\tau$ which occurs
is replaced by its Galois conjugate $-\tau^{-1}$). Thus,
\[ (\tau^{-2})^{\star} = (1 - \tau^{-1})^{\star} = (2 - \tau)^{\star}
= 2 + \tau^{-1} = 1 + \tau = \tau^2. \] 
The corresponding digit matrix, defining the dual substitution
$\zeta^{\star}$, is  
\[ (\D^{\star})^T =  \begin{pmatrix} \{0, \tau^2 \} & 
   \{ 2 \tau^2\} \\ \{ 0 \} & \{ \tau^2 \} \end{pmatrix}. \]
The solution of the corresponding IFS --- with digits 
$\tau^{-2} d^{}_{ijk}$, where $d^{}_{ijk}$ are the digits in 
$(\D^{\star})^T$ --- yields the
prototiles $[0,\tau]$ and $[\tau, \tau+1]$ for $\zeta^{\star}$.
Since this dual tile substitution is nice in the sense that the tiles
are intervals --- rather than some disconnected, fractally shaped sets 
--- it can be written as a word substitution $a^{\star} \mapsto a^{\star}
b^{\star} a^{\star}, \; b^{\star} \mapsto b^{\star} a^{\star}$. 
In particular, even though $\zeta \ne \zeta^{\star}$, they generate
--- up to a scaling by $\tau$ --- the same tiling spaces (the same
tilings), resp. the same biinfinite words.

\section{Dimension 2} \label{secdim2}
In dimension one there is a canonical way to associate a tiling to a
model set, which is, by definition, not a tiling but a uniformly dense
point set:   
The points define a partition of the line into intervals, these
intervals are the tiles. Analogously, a tiling of $\R^1$ by intervals
defines a point set by considering the set of all vertices, i.e.,
the set of all boundary points of the intervals. In 
particular, by assigning each interval to its left boundary point, there 
is a one-to-one correspondence between tiles of the tiling and points 
in the point set. 

In two dimensions the situation becomes more difficult, at least if
one wants to study the different types of tiles/points $x$ through 
their position $x^{\star}$ in the window $W$. It is not obvious which
way to identify tilings with point sets (and vice versa) is the best
one, or the most natural one.  

On one hand, it seems natural to consider the vertices of a tiling
$\T$. In fact, this point of view is taken in many cases, see
e.g. \cite{bksz}, \cite{har}. Then, in general, there is no natural
one-to-one correspondence between tiles  and points. In particular,
the positions $x^{\star}$ in $W$, where $x$ is some vertex, may give
us no relevant information about the tiles in $\T$. One simple reason
is that, in general, the number of tiles per unit area differs 
from the number of vertices per unit area. 

On the other hand, in order to assign a point to each tile, one can use
control points \cite{sol}. Then there is a one-to-one correspondence
between tiles and points. But there are many ways to choose such 
control points, thus one may end up with different windows for different
choices. 

In the sequel we avoid all the problems mentioned above by computing
some  $\star$-dual substitutions of prominent two-dimensional tilings.  
These computations follow the methods in \cite{gel} closely. In the
framework  of Definition \ref{defsubst}, the number of prototiles of a
two-dimensional substitution is considerably
large, since we identify prototiles only up to translation, not up to
isometries. For instance, the Penrose  
tiling uses only two different triangles (up to isometries) as
prototiles, see Figure \ref{fig:penr-subst}. But each of those occurs
in 10 different orientations in 
the tilings. Moreover, we need to distinct a tile and its mirror
image (see the dots in Figure \ref{fig:penr-subst}), thus there are 40
prototiles altogether with respect to translations. Rather than
working with a $40 \times 40$ digit matrix, we will slightly change
our framework: Instead of digits representing translations, we go over
to 'digits' which are isometries. The entries in the according
matrices now are functions rather than translation vectors. Each
function is of the form $x \mapsto QRx + d$, where $Q$ is --- as usual
--- an expanding linear map which is the same for all functions, $R$
is some rotation or reflection, and $d$ is a translation vector, as
above. In other words: We consider {\em matrix function systems}
instead of digit matrices, cf. \cite{lms}. Since we are considering
plane tilings, we will identify $\R^2$ with the complex plane $\C$.
Then, a rotation about the origin is just a multiplication with some
complex number $z$ with $|z|=1$. Each reflection in the complex plane
can be expressed as complex conjugation, followed by some
rotation. The expanding map $Q$ 
will always be a homothety in the remainder of this section.  
Our goal is to express all translations and rotations, and thus all
maps $f^{}_i$ occurring in the IFS under consideration,  in the ring
$\Z[\xi]$, where $\xi$ is some root of unity. To obtain the dual
substitution, we consider the inverse maps $f^{-1}_i$ and apply a 
primitive element of the Galois group $Gal(\Z[\xi]/\Z)$ to them. 
In plain words: we replace each $\xi$ occurring by an appropriate 
Galois conjugate $\xi'$ in order to obtain the desired maps 
$(f^{-1}_i)^{\star}$. This is in fact how the star map acts in the 
following cases, compare \cite{bksz}. To simplify notation, we 
denote $(f^{-1})^{\star}$ by $f^{\sharp}$.

\subsection{The $\star$-dual of the Penrose tiling}
\begin{figure}
\input{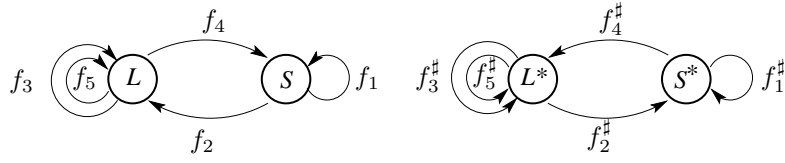}
\caption{The automaton of the IFS of the Penrose tiles (left), compare 
  \eqref{eq:penrifs}, and its dual automaton (right). \label{fig:penrauto}} 
\end{figure}

Consider the prototiles of the Penrose tiling as in Figure
\ref{fig:penr-subst}. It is known that the vertices of a Penrose
tiling are closely related to the cyclotomic field $\Q(\xi^{}_5)$, where
$\xi:= \xi^{}_5 = e^{\frac{2 \pi i}{5}}$. E.g., all vertices, and all edge
lengths, can be expressed in the ring $\Z[\xi]$ of integers in
$\Q(\xi)$. Thus, let the small triangle $S$ have
vertices $0,\; -1$ and $e^{\frac{2 \pi i 3}{10}} = -\xi^4$, and the
large triangle have vertices $0,\; \tau$ and $e^{\frac{2 \pi i}{10}} =
  -\xi^3$. Then, the IFS for these tiles --- arising from the
substitution --- uses the common contracting factor 
$\tau^{-1}=\xi+\xi^4$, and the IFS reads: 
\begin{equation} \label{eq:penrifs}
S = f^{}_1(S) \; \cup \;  f^{}_2(L), \quad L = f^{}_3(L) \; \cup \;
f^{}_4(S) \; \cup \;  f^{}_5(L), 
\end{equation}
where 
\begin{align*} 
f^{}_1(x) & = (\xi+\xi^4) (-\xi^4) x - 1- \xi^3, \\
f^{}_2(x) & =  (\xi+\xi^4) (-\xi) x  -1 - \xi^3 - \xi^4, \\  
f^{}_3(x) & = (\xi+\xi^4) \xi^3 x + 1+ \xi, \\
f^{}_4(x) & =  (\xi+\xi^4) (-\xi^2) \overline{x}  +1 - \xi^3, \\
f^{}_5(x) & =  (\xi+\xi^4) (-1) \overline{x}  +1 - \xi^2 - \xi^3.  
\end{align*}
The corresponding automaton is shown in Figure \ref{fig:penrauto}
(left). Note, that in \cite{gel} the maps are different and do 
{\em not} define the Penrose substitution. Thus, the author of
\cite{gel} did not compute the $\star$-dual (or 'Galois-dual') of the
Penrose substitution, but of a different one, using the same triangles
$S$ and $L$ as prototiles.   

\begin{figure} 
\epsfig{file=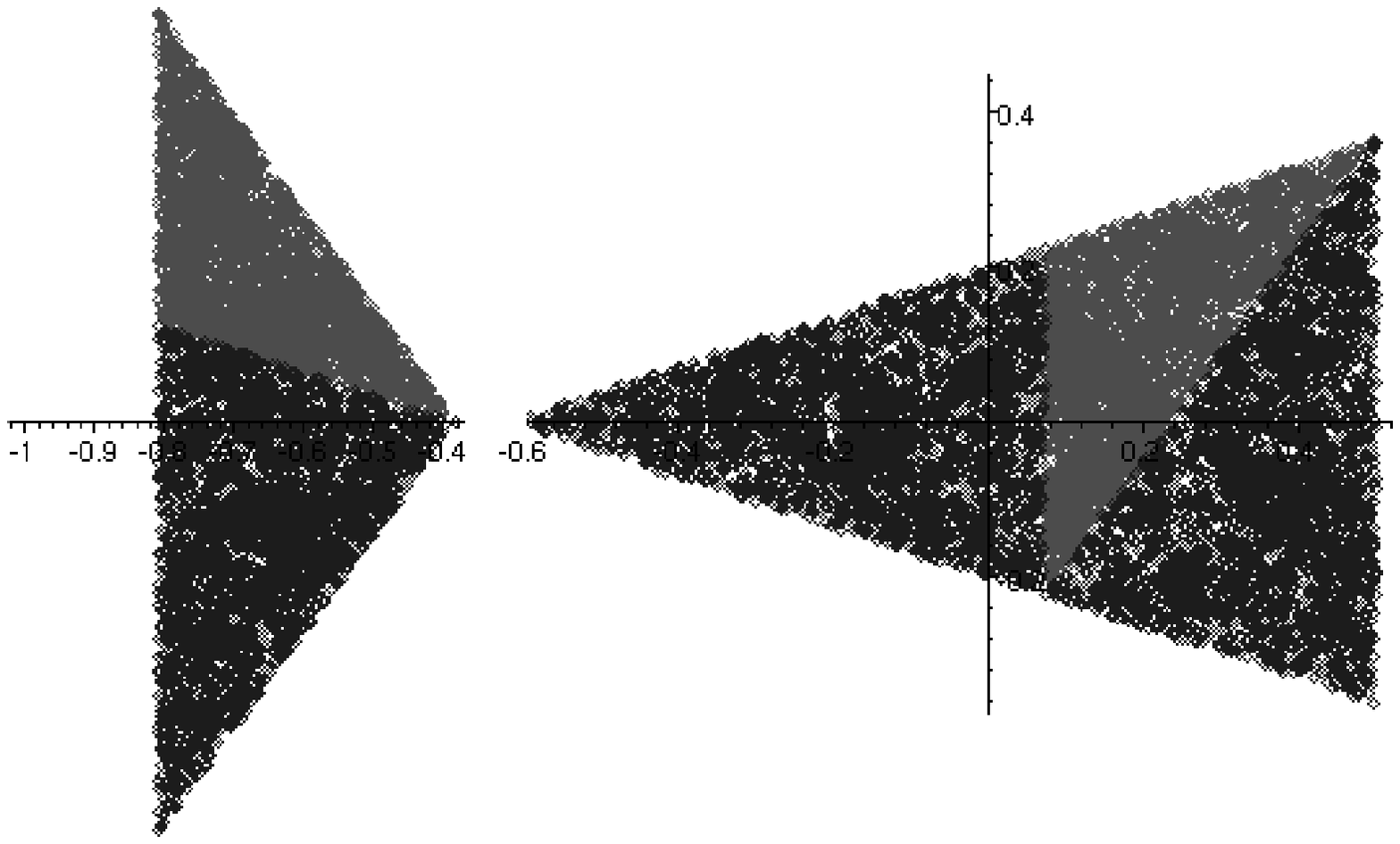, width=70mm} \quad
\epsfig{file=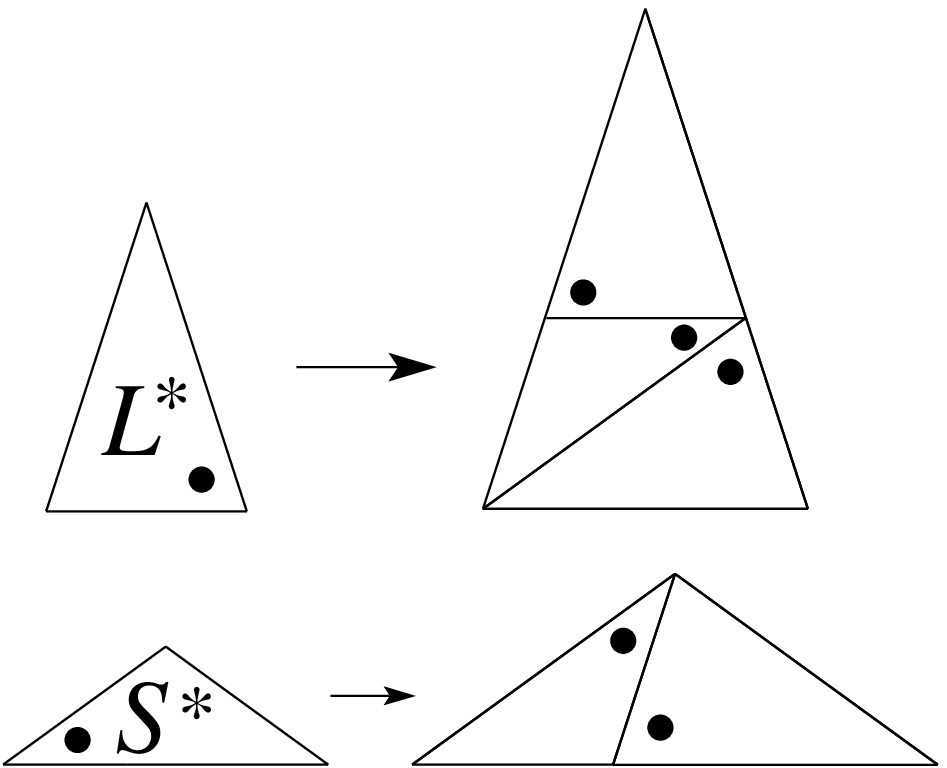, width=50mm}
\caption{A numerical computation of the solution of the dual IFS of
  the Penrose substitution (left), and the T\"ubingen triangle
  substitution (right). One may guess from the images that this is the
  $\star$-dual of the Penrose substitution, which turns out to be
  true. \label{fig:ttt}} 
\end{figure}

In order to calculate the $\star$-dual substitution, consider the 
dual automaton (see Figure \ref{fig:penrauto}, right), resp.\ the dual
IFS. Let $\phi \in Gal(\Z[\xi])$ with $\phi(\xi)=\xi^3$. Note that  
$\phi$ is uniquely defined by this requirement. The $f^{\sharp}_i$ are
obtained by applying $\phi$ to $f^{-1}_i$:
\begin{align*}  
f^{\sharp}_1(x) &= (\xi + \xi^4) (\xi^3 x + \xi^3 +\xi^2) \\
f^{\sharp}_2(x) &= (\xi + \xi^4) (\xi^2x + \xi+\xi^2+\xi^4)\\
f^{\sharp}_3(x) &= (\xi + \xi^4) (-\xi x+\xi+\xi^4) \\
f^{\sharp}_4(x) &= (\xi + \xi^4) (\xi \overline{x} - \xi + \xi^2) \\
f^{\sharp}_5(x) &= (\xi + \xi^4) ( \overline{x} -1+ \xi + \xi^4) \\
\end{align*} 
The corresponding dual IFS reads:
\begin{equation} \label{eq:ifsttt}
S^{\star} = f^{\sharp}_1(S^{\star}) \; \cup \;
f^{\sharp}_4(L^{\star}), \qquad  L^{\star} = f^{\sharp}_2(S^{\star}) \;
\cup \; f^{\sharp}_3(L^{\star}) \; \cup \; f^{\sharp}_5(L^{\star}). 
\end{equation} 
A numerical computation of the solution yields the image in Figure
\ref{fig:ttt} (left). In the figure, the two components of the
solution are torn apart; actually, the solutions do overlap. This fact 
is not relevant, due to Theorem \ref{D+t}. The important point is that
the dissection of each prototile is non-overlapping.
 
It is reasonable to assume that the solution consists of two triangles
similar to the Penrose triangles.  In order to make this precise we
list the coordinates of the triangles: 
\[ L^{\star}: 1+\xi^2+\xi^4, \, 1+\xi+\xi^3, \, 1+\xi^2+\xi^3, \quad
S^{\star}: \xi^2,\xi^3, 2\xi^2 -\xi +\xi^4 -2\xi^3. \] 
All this information provided, it is an easy --- but pretty lengthy --- 
exercise  to check that these triangles are in fact a
solution of \eqref{eq:ifsttt}. The tile-substitution arising from
\eqref{eq:ifsttt} is well-known: It is the T\"ubingen Triangle
tile-substitution (see Figure \ref{fig:ttt}, right). See \cite{fh} for
more details and images of this substitution, as well as several
others. Altogether we have established: 

\begin{thm} \label{penrdual}
The $\star$-dual of the Penrose substitution is the T\"ubingen
triangle substitution, and vice versa. 
\end{thm}

\subsection{The $\star$-dual of the Ammann-Beenker tiling}

\begin{figure}
\epsfig{file=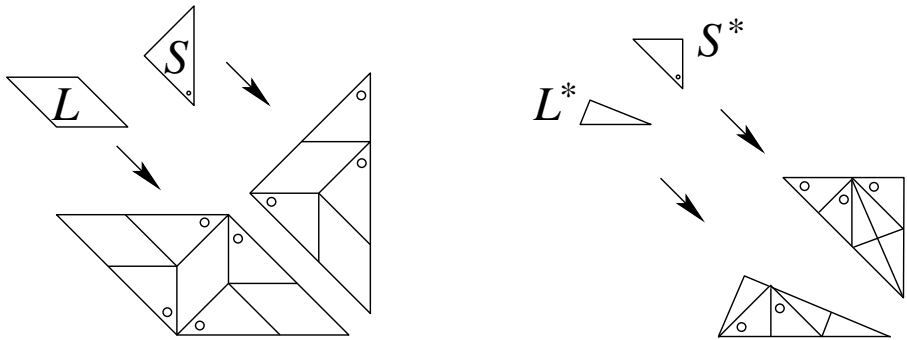, width=120mm}
\caption{The substitution rule for the Ammann-Beenker tiling (left), 
  with prototiles $S,L$, and the star-dual substitution (right), with
  prototiles $S^{\star}, L^{\star}$. Even though the two
  substitutions look different on a first glance, they define very
  similar tilings, see Figure \ref{ab-abd-patch}. 
  \label{ab-abd-subst}} 
\end{figure}

\begin{figure}
\epsfig{file=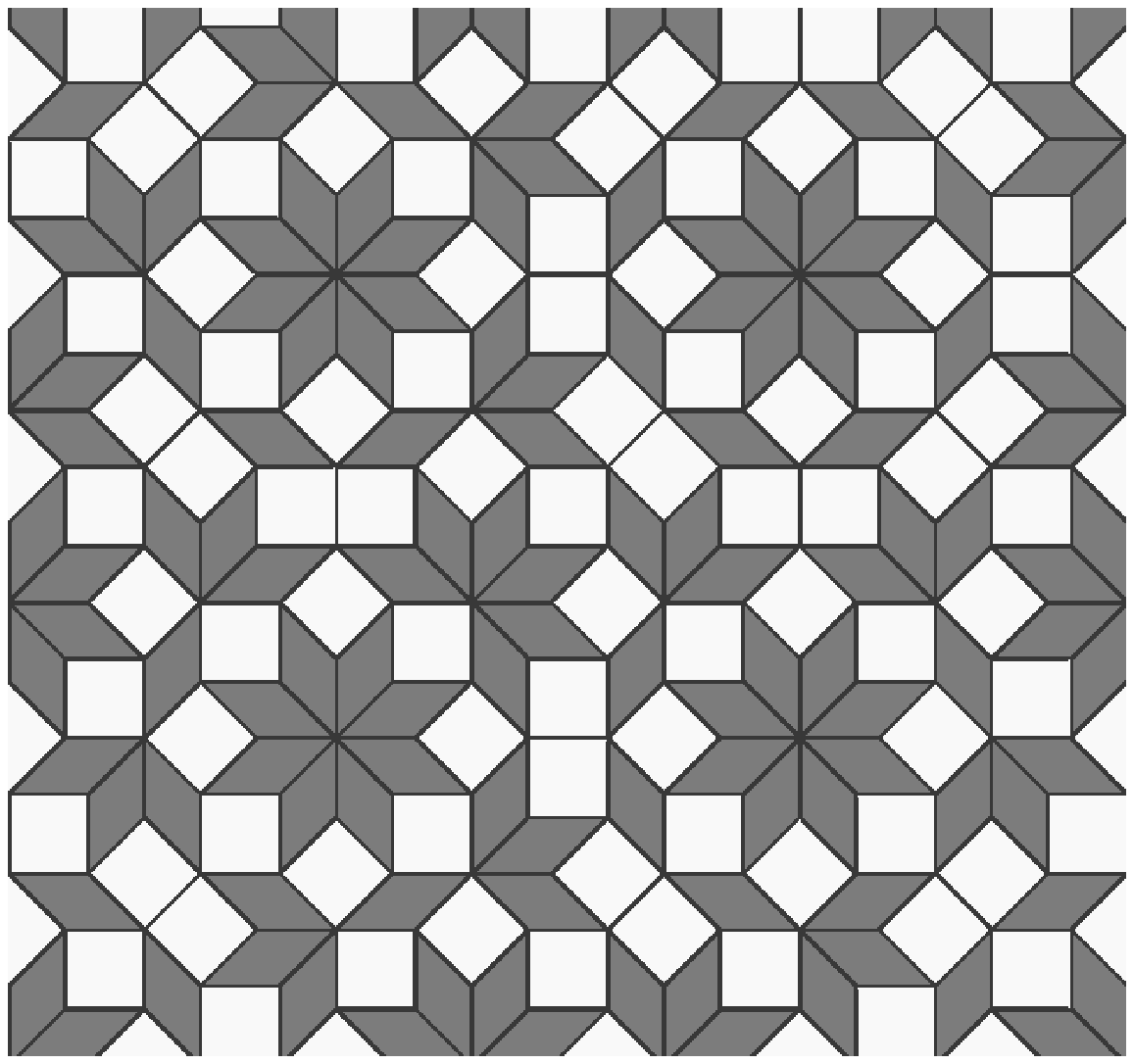, width=60mm}
\epsfig{file=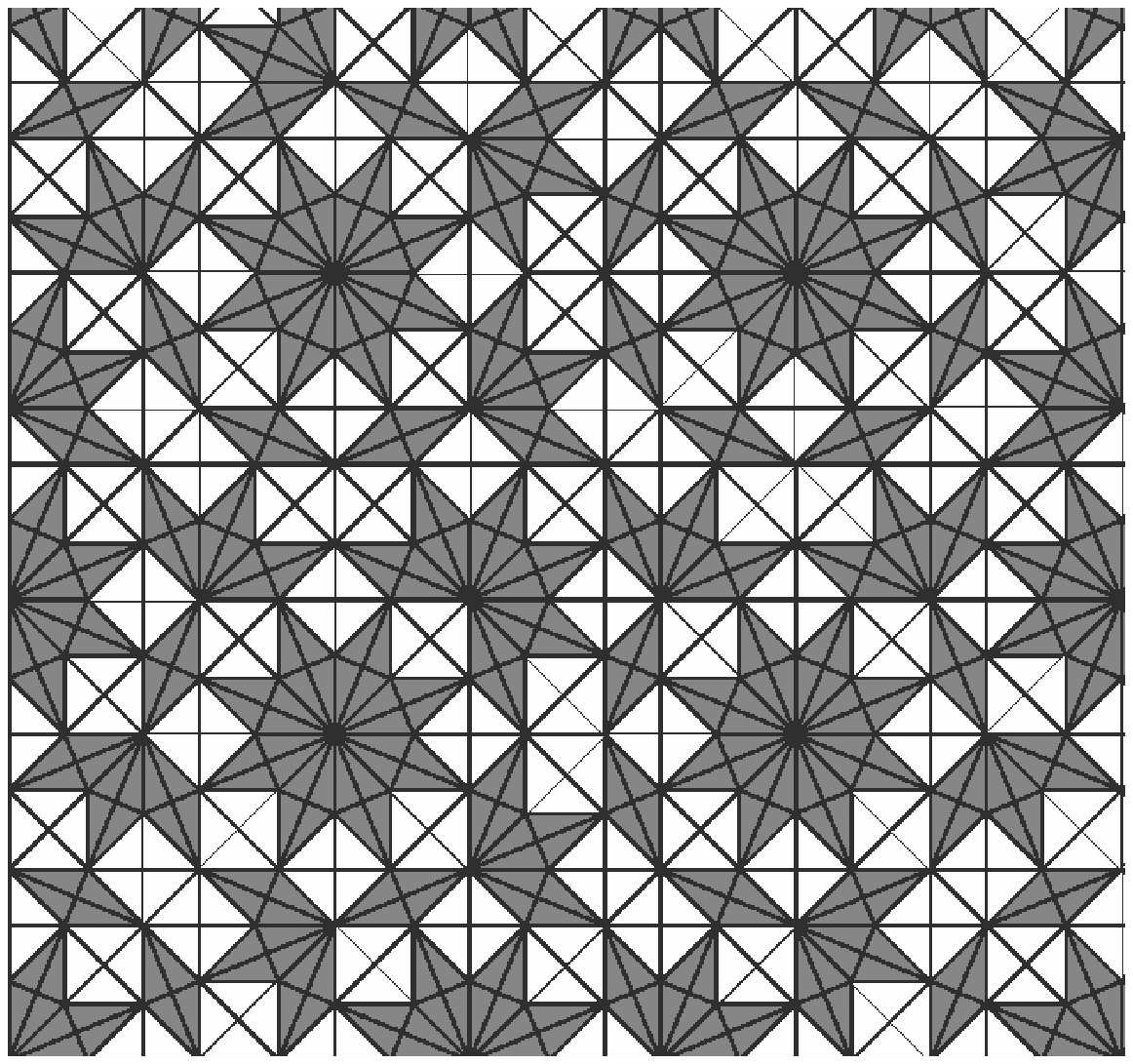, width=60mm}
\caption{Patches of the Ammann-Beenker tiling (left) and its
  star-dual (right). The patches are equivalent in the sense, that one
  obtains the right patch by dividing each tile in the left patch into
  four tiles. \label{ab-abd-patch}}
\end{figure}

In a very similar fashion, one can compute the star-dual of the
Ammann-Beenker tiling \cite{gs}. The substitution rule of the
Ammann-Beenker tiling is shown in Figure \ref{ab-abd-subst} (left). A
part of an Ammann-Beenker tiling is shown in Figure \ref{ab-abd-patch}
(left). The Ammann-Beenker tiling is closely related to the cyclotomic
field $\Q(\xi^{}_8)$, where $\xi : = \xi^{}_8 = e^{\frac{2 \pi i}{8}}$. 
Its $\star$-dual can be calculated in a very similar fashion as in the
case of the Penrose tiling. 
The substitution factor is $\lambda = 1 + \sqrt{2} = 1 +
\xi + \xi^7 = 1 + \xi + \xi^{-1}$, its inverse $\lambda^{-1} = -1 +
\xi +\xi^{-1}$. Let $S$ be the triangular prototile, with vertices $0, 
i, 1+i$, and let $L$ the rhombic prototile with vertices $0,1,
1+\xi,\xi$. The IFS for these prototiles is 
\begin{align*}
L & = f^{}_1(L) \cup f^{}_2(L) \cup f^{}_3(L) \cup f^{}_6(S) \cup
f^{}_{7}(S) \cup f^{}_{8}(S) \cup f^{}_{9}(S), \\
S & =  f^{}_4(L) \cup f^{}_5(L) \cup f^{}_{10}(S) \cup f^{}_{11}(S)
\cup f^{}_{12}(S),  
\end{align*} 
where 
\begin{align*}
f^{}_1(x) & = \lambda^{-1} \xi^4 x +1 +\xi, & 
f^{}_7(x) & = \lambda^{-1} \xi^3 x + 1 + 2\xi + \xi^2 + \xi^7, \\
f^{}_2(x) & = \lambda^{-1} \xi   \overline{x} +1 +\xi +\xi^2 +\xi^7,& 
f^{}_8(x) & = \lambda^{-1} \xi^6 \overline{x}+2 +\xi +\xi^2 +\xi^7,\\
f^{}_3(x) & = \lambda^{-1} \xi^7 \overline{x} +1 +\xi +\xi^2, & 
f^{}_9(x) & = \lambda^{-1} \xi^7 x  +1, \\
f^{}_4(x) & = \lambda^{-1}       x +\xi +\xi^2, & 
f^{}_{10}(x) & = \lambda^{-1} \xi^5 x + \xi + \xi^3, \\
f^{}_5(x) & = \lambda^{-1} \xi^2 x +\xi, &
f^{}_{11}(x) & =  \lambda^{-1} \xi^6 \overline{x} + 1 + \xi + \xi^2, \\
f^{}_6(x) & = \lambda^{-1} \xi^2 \overline{x} +\xi, & 
f^{}_{12}(x) & =  \lambda^{-1} \xi^3 x + 2 \xi + \xi^2. \\
\end{align*}

The functions $f^{\sharp}_i$ of the dual IFS are easily obtained by
considering the inverse functions $f^{-1}_i$, and replacing each $\xi$
occurring  by its Galois conjugate $\xi^3$ (or $\xi^5$ or $\xi^7$, the
resulting dual tiling is always the same). The dual IFS reads:
\begin{align*}
L^{\star} & = f_1^{\sharp}(L^{\star}) \cup f_2^{\sharp}(L^{\star})
\cup f_3^{\sharp}(L^{\star}) \cup f_4^{\sharp}(S^{\star}) \cup
f_5^{\sharp}(S^{\star}), \\
S^{\star} & = f_6^{\sharp}(L^{\star}) \cup f_7^{\sharp}(L^{\star})
\cup f_8^{\sharp}(L^{\star}) \cup f_9^{\sharp}(L^{\star}) \cup
f_{10}^{\sharp}(S^{\star}) \cup f_{11}^{\sharp}(S^{\star}) \cup
f_{12}^{\sharp}(S^{\star}). 
\end{align*}
Its solution are the two tiles shown in Figure \ref{ab-abd-subst}
(right): $S^{\star}$ is an isosceles orthogonal triangle, and 
$L^{\star}$ is a thin orthogonal triangle. Note that the original
triangular prototile $S$ can be dissected into two copies of
$S^{\star}$, and the rhombic prototile $L$ can be dissected into four
copies of $L^{\star}$. In fact, unlike in the Penrose case, the 
$\star$-dual tiling of the Ammann-Beenker tiling is closely related to
the Ammann-Beenker tiling itself, compare Figure \ref{ab-abd-patch}: 
dissecting each $L$ into four tiles, and each $S$ into two tiles, 
yields the $\star$-dual tiling of the Ammann-Beenker
tiling. Altogether we have established:

\begin{thm} \label{ammbdual}
The $\star$-dual of the Ammann-Beenker tiling is mutually
locally derivable with the Ammann-Beenker tiling. 
\end{thm}

\section{Self-duality} \label{secselfd}

The last example, the Ammann-Beenker tiling and its $\star$-dual,
motivates the question whether there are substitutions which are
strictly self-dual.  

\begin{defi}
Let $\sigma$ be a tile-substitution such that $\X^{}_{\sigma}$ consists 
of cut-and-project tilings. If $\X^{}_{\sigma} =
\X^{}_{\sigma^{\star}}$, then $\sigma$ is called {\em self-dual (with
respect to $\star$-duality)}. 
\end{defi}

Thus, the example in Section \ref{secstardual}, Fibonacci squared, is
self-dual.  
Since the star-dual tiling of a tiling in $G$ is a tiling in $H$ (cf.\
Definition \ref{defmodset}), a trivial necessary condition for a tiling to
be self-dual is $G=H$. Hence, if the substitution factor $\lambda$ of a
self-dual substitution is a unimodular PV-number, then $\lambda$ must
be a quadratic algebraic number. A refined necessary condition is
obtained by considering the substitution matrix $M^{}_{\sigma}$. 

\begin{lem} \label{lemsd}
If $\sigma$ is a self-dual substitution in $\R^d$, then the following 
property holds for the substitution matrix $M^{}_{\sigma}=(|\D^{}_{ij}|)$:
\begin{equation} \label{pmp}
  (M^{}_{\sigma})^T = P  M^{}_{\sigma} P^{-1}, 
\end{equation}
where $P$ is some permutation matrix.
\end{lem}
\begin{proof} 
The result follows immediately from $M^{}_{\sigma^{\star}} =
\big( |(\D^{}_{ij})^{\star}| \big)^T$, by definition of
$\sigma^{\star}$. The tiling spaces $\X_{\sigma}$ defined by a
substitution are unique up to enumeration of the tiles, thus the
permutations enter. 
\end{proof}

In fact, permutations has to be taken into account.
There are self-dual substitutions where $(M^{}_{\sigma})^T \ne
M^{}_{\sigma}$; in other words: where $M^{}_{\sigma}$ is not symmetric.
For instance, the dual substitution of $\sigma: \; a \to ab, \; b \to
aab$ is $\sigma^{\star}: \; a^{\star} \to b^{\star} b^{\star}
a^{\star}, \; b^{\star} \to b^{\star} a^{\star}$. By identifying the
first tile --- $a$ --- of $\sigma$ with the second tile --- $b^{\star}$
--- of the dual substitution $\sigma^{\star}$, one obtains that the two
substitutions define the same tiling spaces. Hence $\sigma$ is
self-dual. The substitution matrix
$M^{}_{\sigma} = \big( \begin{smallmatrix} 1 & 1 \\ 2 & 1
\end{smallmatrix} \big)$ is not symmetric, but 
\[ P^{-1} M^{}_{\sigma} P =  \begin{pmatrix} 0 & 1 \\ 1 & 0 \end{pmatrix}
  \begin{pmatrix} 1 & 1 \\ 2 & 1 \end{pmatrix}  \begin{pmatrix} 0 & 1 \\
  1 & 0 \end{pmatrix} =  \begin{pmatrix} 1 & 2 \\ 1 & 1 \end{pmatrix}
  = (M^{}_{\sigma})^T, 
\]
thus \eqref{pmp} holds. 

A sufficient condition for a substitution to be self-dual is an
immediate consequence of Theorem \ref{D+t}. 

\begin{cor} \label{corsd}
Let $\sigma$ be a substitution with digit matrix $\D$ and
$\sigma^{\star}$ its dual substitution, with digit matrix $\D' =
(\D^{\star})^T$. If there are vectors $t^{}_1,t^{}_2, \ldots, t^{}_m
\in \R^d$, $A := (t^{}_j - Q t^{}_i)_{1 \le i,j  \le  m}$ and a
permutation matrix $P$ such that 
\[ P \D' P^{-1} = \D + A , \]
then $\sigma$ is self-dual. \hfill $\square$
\end{cor}

This result gives a purely algebraic test for self-duality.

\section{Outlook}
The notion of star-duality is a concept tailored to arbitrary model
sets in any locally compact Abelian group. Restricted to the
unimodular case, it coincides 
with other notions of duality, in particular the 'natural
decomposition method' in \cite{sirw}, and the dual tilings considered
by P.\ Arnoux, S.\ Ito and others, see for instance
\cite{ei}, \cite{sai}. This is the reason why we were lazy at some
occasions, using the term 'self-dual' instead of 'self-$\star$-dual' 
or 'self-dual w.r.t.\ star-duality'. Theorem \ref{D+t} can be utilised
to show the equivalence of the concepts mentioned. 

The benefit of star-duality seems to be its generality, and its
accessibility by purely algebraic methods. This paper tries to make
first steps in this direction.  
For future work, it should be promising to apply star-duality to the case
of word substitutions on two or three letters. In the case of Sturmian
substitutions (symbolic substitutions on two letters, where the
corresponding biinfinite words contain exactly $n+1$ subwords of
length $n$, for all $n\ge 0$, see \cite{fogg}) the results in
\cite{blr} can be used to give a complete characterisation of
self-dual (Sturmian) substitutions. For details we refer to future
work.  

A further question is: are there self-dual substitutions in dimension
$d \ge 2$? The Ammann-Beenker tiling is very close to be self-dual, but
is is not exactly self-dual. In dimension one, there are several
self-dual substitutions. Trivially, these examples generalise to higher 
dimensions, by using the Cartesian product (for instance, Fibonacci
times Fibonacci). But it remains to find a two-dimensional self-dual
substitution apart from that. 

\section*{Acknowledgements}

It is a pleasure to thank V.\ Sirvent for cooperation and helpful
comments. The author is indebted to the referee for many 
valuable remarks. This work was supported by the 
German Research Council (DFG), within the CRC 701.

\end{document}